\documentclass{article}
\author{Romain Bondil}
\title{Discriminant of a generic projection of  a minimal normal surface singularity} 

\usepackage{amsmath,amssymb,amsthm,amsfonts,latexsym}
\usepackage[english]{babel}
\input utile.def
\newcommand{\C}{\mathbb{C}}
\def\DebEq{\vskip-8pt}
\def\FinEq{\vskip-3pt \noindent}
\newtheorem{theorem}{Theorem}[section]
\newtheorem{lemma}[theorem]{Lemma}
\newtheorem{e-proposition}[theorem]{Proposition}
\newtheorem{corollary}[theorem]{Corollary}
\newtheorem{e-definition}[theorem]{Definition\rm}
\newtheorem{notation}[theorem]{Notation\rm}
\newtheorem{remark}{\it Remark\/}
\newtheorem{example}{\it Example\/}

\newtheorem{f-notation}[theoreme]{Notation}

\begin{document}

\author{Romain Bondil}
\date{February 2003}

\maketitle
\begin{abstract}
Let $(S,0)$ be a rational complex  surface singularity with reduced fundamental cycle, also known as  a {\em minimal} singularity. Using a fundamental result by M. Spivakovsky, we explain how to get a minimal resolution of the discriminant curve for a generic projection of $(S,0)$ onto $(\C^2,0)$ from the resolution of $(S,0)$.
\end{abstract}

The material in this Note is organized as follows~: \S~\ref{s-discri} recalls the definitions of polar curves, discriminants and a remarkable property of transversality due to Brian\c{c}on-Henry and Teissier (thm.~~\ref{thm-transv-polaire}). For minimal surface singularities, a theorem due to M.~Spivakovsky describes the behavior of the generic polar curve (cf. \S~\ref{s-mini}). We use this theorem in  \S~\ref{s-key-lemmas} to prove two lemmas relating on the one side the resolution of the generic polar curve to the resolution of a minimal surface singularity, and on the other side, the polar curve and the discriminant. Gathering these results, we give our main theorem  in \S~\ref{s-main-result}, which provides us with a combinatorial way to describe  the discriminant.

\section{Polar curves and discriminants}
\label{s-discri}
Let $(S,0)$ be a normal complex surface singularity $(S,0)$, embedded in $(\C^N,0)$: for any $(N-2)$-dimensional vector subspace $D$  of $\C^N$, we consider a linear projection $\C^N \rightarrow \C^2$ with kernel $D$ and denote by 
$p_D : (S,0)\rightarrow (\C^2,0),$ the restriction of this projection to $(S,0)$. 

Restricting ourselves to the $D$ such that $p_D$ is finite, and considering a small representative $S$ of the germ $(S,0)$, we define, as in \cite{pol-loc}~(2.2.2), the {\em polar curve} $C_1(D)$ of the germ $(S,0)$ relative to the direction $D$, as the closure in  $S$ of the critical locus of the restriction of $p_D$ to $S\setminus\{0\}$. As explained in loc. cit., it makes sense to say that for an open dense subset of the Grassmann manifold $G(N-2,N)$  of  $(N-2)$-planes in $\C^N$, the space curve $C_1(D)$ are {\em equisingular} in term of strong simultaneous resolutions.

Then we define the {\em discriminant} $\Delta_{p_D}$ as (the germ at $0$ of) the reduced analytic curve of $(\C^2,0)$ image of $C_1(D)$ by the finite morphism $p_D$. 

Again, one may show that,  for a generic choice of $D$, the discriminants obtained are {\em equisingular germs of plane curves}, but we will need a much more precise result, that demands  another definition (cf. \cite{BGG}~IV.3)~: 

\begin{e-definition}
\label{def-generic-proj-curve}
{\rm Let $(X,0)\subset (\C^N,0)$  be a germ of reduced curve. Then a linear projection $p : \C^N \rightarrow \C^2$ will said to be {\em generic}  w.r.t. $(X,0)$ if the kernel of $p$ does not contain any limit of secants to $X$ (cf. \cite{BGG} for an explicit description of  the cone $C_5(X,0)$ formed by the limits of secants to $(X,0)$).}
\end{e-definition}

We now  state the following transversality result (proved for curves on surfaces of $\C^3$ in \cite{BH} thm.~3.12  and in general as the ``lemme-cl\'e'' in \cite{laRabida}~V~(1.2.2))~:

\begin{theorem}
\label{thm-transv-polaire}
Let $p_D~: (S,0)\rightarrow (\C^2,0)$ be as above, and $C_1(D)\subset (S,0)\subset (\C^N,0)$ be the corresponding polar curve, then there is an open dense  subset $U$ of $G(N-2,N)$ such that for $D\in U$ the restriction of $p$ to $C_1(D)$ is generic in the sense of def.~\ref{def-generic-proj-curve}.
\end{theorem}

\begin{e-definition}
\label{def-gen-discri}{\rm 
For all $D$ in the open subset $U$ of thm.~\ref{thm-transv-polaire},  the discriminant $\Delta_{p_D}$ are  equisingular in the sense of the well-known equisingularity theory for germs of plane curves (cf. e.g. the account at the beginning of~\cite{BGG})~: we will call this equisingularity class  the {\em generic discriminant} of $(S,0)$.}
\end{e-definition}

\section{Polar curves for minimal singularities of surface after Spivakovsky}
\label{s-mini}


We first recall how one may define a minimal singularity in the case of normal surfaces~ (cf.~\cite{mark} II.2)~:

\begin{e-definition}
\label{def-mini}{\rm 
A normal surface singularity $(S,0)$ is said to be {\em minimal} if it is rational with reduced fundamental cycle (see \cite{Artin} for these latter notions).} 
\end{e-definition}

Let $\pi~: (X,E) \rightarrow (S,0)$ be the minimal resolution of the singularity $(S,0)$, where $E=\pi^{-1}(0)$ is the exceptional divisor, with components $L_i$. A {\em cycle} will be by definition a divisor with support on $E$ i.e. a linear combination $\sum a_i L_i$ with $a_i\in \ZZ$ (or $a_i\in \QQ$ for a $\QQ$-cycle).


Considering the dual graph $\Gamma$ associated to the exceptional divisor $E$ (cf. \cite{mark} I.~1) in which each component $L_x$ gives a vertex $x$ and two vertices are connected if, and only if, the corresponding components intersect, the minimal singularities have the following easy characterization (cf. loc. cit. II.~2)~:
\begin{lemma}
\label{lem-comb-char}
Let $\Gamma$ be as above the dual graph associated to the minimal resolution of a normal surface singularity $(S,0)$. For  each vertex $x\in \Gamma$, one defines its weight $w(x):=-L_x^2$ (self-intersection of the corresponding component $L_x$) and its valence $\gamma(x)$ which is the number of vertices connected to $x$.
Then $(S,0)$ is minimal if, and only if, $\Gamma$ is a tree and for all $x\in \Gamma, \: w(x)\geq \gamma(x)$. 
\end{lemma}

To two vertices $x,y\in \Gamma$ (which is a tree), we associate the shortest chain in $\Gamma$ connecting them, which we denote by $[x,y]$. The {\em distance} $d(x,y)$ is by definition the number of edges on $[x,y]$.

In~\cite{mark}~III.~5, generalizing an earlier work by G. Gonzalez-Sprinberg in \cite{nash-pdr}, M. Spivakovsky further introduces the following number $s_x$ associated to each vertex $x\in \Gamma$.
If $Z.L_x<0$ (where $.$ denotes the intersection number) then put $s_x:=1$ (and $x$ is said to be non-Tyurina). Otherwise $x$ is said to be a Tyurina vertex, then denote $\Delta$ the Tyurina component of $\Gamma$ containing $x$ (i.e. the maximal connected subgraph of $\Gamma$ containing only Tyurina vertices), and put $s_x:=d(x,\Gamma\setminus\Delta)+1$.

These numbers $s_x$ coincide, in the special case of minimal singularities, with the  {\em desingularization depths} introduced in \cite{Le-meral} p.~8.

Let $x,y$ be two adjacent vertices~: the edge ($x,y)$ in $\Gamma$ is called a {\em central arc} if $s_x=s_y$. A vertex $x$ is called a {\em central vertex} if there are at least two vertices $y$ adjacent to $x$ such that $s_y=s_x-1$ (cf. loc. cit.).

Eventually, we define the following $\QQ$-cycle $Z_\Omega$ on the minimal resolution $X$ of $(S,0)$ by~:
\DebEq
\begin{equation}
\label{eq-ZOmega}
Z_\Omega=\sum_{x\in \Gamma} s_x L_x - Z_K,
\end{equation}\FinEq
\noindent where $\Gamma$ is the dual graph of the resolution, and $Z_K$ is the numerically canonical $\QQ$-cycle \footnote{uniquely defined by the condition that for all $x\in \Gamma,  Z_K.L_x=-2-L_x^2$ since the intersection product on $E$ is negative-definite}.

One may now quote the important theorem 5.4 in loc. cit. in the following way~:\footnote{see also the account in~\cite{Le-sand} (7.4), just beware that one term is missing in the formula giving $m_x:=-Z_\Omega.L_x$ there.}
\begin{theorem}
\label{thm-mark}
Let $(S,0)$ be a {\em minimal} normal surface singularity. There is a open dense subset $U'$ of the open set $U$ of thm.~\ref{thm-transv-polaire}, such that for all $D\in U'$ the strict transform $C'_1(D)$ of $C_1(D)$ on $X$~:

\noindent a) is a multi-germ of smooth curves intersecting each component $L_x$ of $E$ transversally in exactly $-Z_\Omega.L_x$ points,

\noindent b) goes through the  point of intersection of $L_x$ and $L_y$ if and only if $s_x=s_y$ (point corresponding to a {\em central arc} of the graph). Further, the $C_1'(D)$, with $D\in U'$  do not share other common points (base points) and these base points   are simple,  i.e. the $C_1'(D)$ are separated when one blows-up these points once.
\end{theorem}

\section{From polar curves to discriminants, key lemmas}
\label{s-key-lemmas}
\begin{lemma}
\label{lem-res-cp}
Let $(S,0)$ be a {\em minimal} normal surface singularity, embedded in $\C^N$ and $\pi~: X\rightarrow (S,0)$ its minimal resolution. It is known that $\pi$ is (the restriction to $S$ of) a composition  $\pi_1\circ\cdots\circ \pi_r$ of point blow-ups in~$\C^N$. We claim that this composition of blow-ups  is also the  minimal resolution of the generic polar curve $C_1(D)$ for $D\in U'$ as in thm.~\ref{thm-mark}.
\end{lemma}

\begin{proof}
  The fact that $\pi$ is a composition of point blow-ups  is general for rational surface singularities (for a non-cohomological proof in the case of minimal singularities, see \cite{RES}~5.9). Conclusion a) in thm.~\ref{thm-mark} certainly gives that $\pi$ is a resolution of $C_1(D)$. We prove that this resolution is minimal~: among the exceptional components in $X$ obtained by the last point blow-up, there is  a component $L_x$ corresponding either to a central vertex of $\Gamma$ or to the boundary of a central arc.

If $L_x$ corresponds to a central vertex, one computes from (\ref{eq-ZOmega}) page~\pageref{eq-ZOmega}, the number of  branches of  $C_1'(D)$ intersecting $L_x$, i.e. $-Z_\Omega.L_x=-( \sum s_y + (s_x+1) L_x^2 +2)$. By the definition of a central vertex (before thm.~\ref{thm-mark}), this  must be at least two, which proves that these  branches are not separated before $L_x$ is obtained.

If  $L_x$ is the boundary of a central arc, let $L_y$ be the other boundary~: then both $L_x$ and $L_y$ appear as exceptional components of the last blow-up $\pi_r : X\rightarrow S_{r-1}$ at $0_{r-1}$.  Now, the strict transform of $C_1(D)$ at the point $0_{r-1}$ can not be smooth. Indeed,  by an argument in \cite{LJ-GS}~1.1, if it were smooth, then its strict transform $C_1'(D)$ on $X$, smooth surface,  would go through a smooth point of the exceptional divisor.
\end{proof}

\begin{lemma}
\label{lem-mult-two}
For $D\in U'$ as in thm.~\ref{thm-mark}, the polar curve $C_1(D)$ on ($S,0)$ has only  smooth branches and branches of multiplicity two, the latter being exactly those whose strict transform go through a central arc as in b) of theorem~\ref{thm-mark}.
\end{lemma}

\begin{proof}
For $D\in U'$ one may compute the multiplicity $e(\Delta_{p_D},0)$ of the discriminant en $0$ by the following (cf. e.g. \cite{RES} \S~4.4~: here the divisorial discriminant is reduced by genericity)~:\DebEq
\begin{equation}
\label{eq-mult-discri}
e(\Delta_{p_D},0)=\mu-1+e(S,0),
\end{equation}\FinEq
where $\mu$ is the Milnor number of a generic hyperplane section of $(S,0)$ and $e(S,0)$ is the multiplicity of the surface.
A well-known formula (see e.g. \cite{mynote} prop.~5, taking $I=m$) for $\mu$ reduces, for $(S,0)$ minimal, to $\mu=1+Z.Z_K$. In turn, 
$e(S,0)=-Z^2$ (cf. \cite{Artin}), which in~(\ref{eq-mult-discri}) reads~:
\DebEq
\begin{equation}
\label{eq-mult-comb}
e(\Delta_{p_D},0)=Z.Z_K-Z^2.
\end{equation}
 \FinEq
On the other hand, the number $n_b$ of branches of $\Delta_{p_D}$ is described by thm.~\ref{thm-mark}~: denote by $n_{bs}$ the number of those branches which go through a central arc and so are counted twice in $Z.Z_\Omega$, then~:
\DebEq
\begin{equation}
\label{eq-nb}
n_b=-Z.Z_\Omega-n_{bs}.
\end{equation}
\FinEq
Using expression~(\ref{eq-ZOmega}) for $Z_\Omega$ in (\ref{eq-nb}), and then using (\ref{eq-mult-comb}) yields~:
$n_b=-\sum_{x\in \Gamma} s_x L_x.Z +Z.Z_K -n_{bs}$ and
$n_b=e(\Delta,0)-n_{bs}- \sum_{x\in \Gamma} (s_x-1) Z.L_x$,
but, by definition $s_x=1$ if $Z.L_x\neq 0$, whence~:
\DebEq
\begin{equation}
\label{eq-final-nb}
n_b=e(\Delta,0)-n_{bs}.
\end{equation}
\FinEq
Since we know from the proof of lemma~\ref{lem-res-cp} that all the branches of $C_1(D)$ counted in $n_{bs}$ are actually singular, (\ref{eq-final-nb}) proves the whole assertion of the current lemma.
\end{proof}

\begin{corollary}
\label{cor-res-discri}
Take the chain of point blow-ups over $(\C^N,0)$ that gives the minimal resolution of $(C_1(D),0)$ for $D\in U'$. Then performing over $(\C^2,0)$ the ``same'' succession of blow-ups (this makes sense because of footnote~\ref{ft-contact}), we get the minimal resolution of the plane curve $\Delta_{p_D}=p_D(C_1(D))$.
\end{corollary}
\begin{proof}
Since, by lem.~\ref{lem-mult-two}, the multiplicity of the branches of $C_1(D)$ is at most two, these branches are plane curves and so are equisingular to their {\em generic} projection by $p_D$ (here we use thm.~\ref{thm-transv-polaire})~: so much for the branches. Further, by another result of Teissier's (see \cite{laRabida} Chap.~I (6.2.1) and remark p. 354) a generic projection is bi-lipschitz, which implies that it preserves the contact between branches.\footnote{ \label{ft-contact}Indeed, the contact between two branches $\gamma_1(t)$ and $\gamma_2(t)$ which are both of multiplicity one or two, that we define as the number of blow-ups to separate them, may be read from the order in $t$ of the difference $\gamma_1(t)-\gamma_2(t)$, which is a bi-lipschitz invariant.
 Since we blow-up always in the ``same chart'' these blow-ups  actually dominate the blow-ups of the plane, as claimed in the corollary.}
\end{proof}

\section{Statement of the main result}
\label{s-main-result}
Gathering the results from lemma~\ref{lem-res-cp} to  cor.~\ref{cor-res-discri}, we obtain~:
\begin{theorem}
\label{thm-surf-disci}
Let $(S,0)$ be a minimal normal surface singularity, embedded in $\C^N$ and $\pi~: X\rightarrow (S,0)$ its minimal resolution, which is a composition  $\pi_1\circ\cdots\circ \pi_r$ of point blow-ups in $\C^N$. Let $\Delta_{S,0}$ be the generic discriminant of $(S,0)$ (cf. def.~\ref{def-gen-discri}). Then, performing over $(\C^2,0)$ the ``same'' succession of blow-ups (cf. footnote~\ref{ft-contact}), we get the minimal resolution of the plane curve $\Delta_{S,0}$.
\end{theorem}

We claim that this result, together with thm.~\ref{thm-mark}  gives an easy way to get a combinatorial description of $(\Delta_{S,0},0)$~:

\begin{notation}
\label{not-delta}
{\rm 
i) We denote by $\Delta_{A_n}$ the generic discriminant of the $A_n$ surface singularity, which is the equisingularity class of the plane curve defined by $x^2+y^{n+1}=0$.

ii) We denote by $\delta_n$ the generic discriminant of the singularity which is a cone over a rational normal curve of degree $n$ in $\PP^n_\C$~: it is defined by $2n-2$ distinct lines through the origin.}
\end{notation}
\smallskip

The assertion in ii) follows from the fact that $C_1(D)$ is the cone over the critical locus of the projection from the rational normal curve onto a line, which has degree $2n-2$ by Hurwitz formula.

We need to introduce several subsets of a dual graph $\Gamma$~: we denote by $\Gamma_{NT}=\{x_1,\dots,x_n\}$ the set of Non-Tyurina vertices in  $\Gamma$, which are here the $x\in \Gamma$ such that $w(x)>\gamma(x)$ (notation as in lem.~\ref{lem-comb-char}).

We denote by $\cC_v$ resp. $\cC_a$ the set of central vertices resp. central arcs in $\Gamma$ (cf. def. before thm.~\ref{thm-mark}).

\begin{corollary}
From thm.~\ref{thm-mark} we know that the components of the strict transform $C_1(D)'$  of $C_1(D)$ on the resolution $X$ of $(S,0)$ go through components corresponding to elements  of $\Gamma_{NT}\cup \cC_a\cup \cC_v$, and we also know the number of branches of $C_1(D)'$ on each of these components.

From thm.~\ref{thm-surf-disci} we know the contact between the corresponding branches of $C_1(D)$ (or $\Delta_{S,0}$)~: the contact between two branches whose strict transforms lie  respectively on a component $L_x$ and a component $L_y$ equals $1+N$,  where $N$ is the number of blow-ups  necessary so that $L_x$ and $L_y$ no longer be in the same Tyurina component of the corresponding $(S_N,0_N)$ singularity, with the further requirement that if ,say, the first branch actually goes through  a central arc $L_x\cap L_{x'}$,  the number $N$ corresponds to the number of blow-ups so that {\em both} $x$ and $x'$ no longer be in the same Tyurina component as $y$.
\end{corollary}

From this, we easily see that each $x_i\in \Gamma_{NT}$ contributes with a $\delta_{x_i}:=\delta_{w(x_i)-\gamma(x_i)}$ (cf.~\ref{not-delta} ii)), i.e. $2(w(x_i)-\gamma(x_i))-2$ lines, and that the contact between these $\delta_{x_i}$ and other branches of the discriminant is one.
For the contribution of the central elements, we first compute the number of branches on each components with thm.~\ref{thm-mark} and use thm.~\ref{thm-surf-disci} for the contact as in the following examples~:

\begin{example}
\label{ex-graphe}{\rm
Consider $(S,0)$ with the graph $\Gamma$ as on figure~1 below, where the $\bullet$ denote Tyurina vertices (with $w(x)=\gamma(x)$), and $\Gamma_{NT}=\{x_1,\dots,x_4\}$ with the weights indicated on the graph.  Remark  that as a general rule $\delta_{x_i}=\emptyset$ when $w(x_i)=\gamma(x_i)+1$, hence here only $x_1$ actually gives a $\delta_{x_1}$equals to four lines.

\noindent i) In the {\em first Tyurina component} (bounded by $x_1,x_2,x_4$)
there is a central vertex and a central arc, which respectively give a $\Delta_{A_5}$ and a $\Delta_{A_4}$ curve.


 After two blow-ups the boundaries of the central arc and the central vertex are in distinct Tyurina components, hence the contact between the $\Delta_{A_5}$ and $\Delta_{A_4}$ is three.

\noindent ii) In the {\em second Tyurina component} (bounded by $x_2$ $x_3$), there is a central vertex : this gives a $\Delta_{A_3}$ which has contact $1$ with the others $\Delta_{A_i}$ obtained.

Hence, using coordinates, we may give as representative of  the equisingularity class of $\Delta_{S,0}$~:

$\Delta_{S,0}\: : \: {(x^4+y^4)}(x^2+y^6)(x^2+y^5)(y^2+x^4)=0.$
}
\end{example}

\medskip

\setlength{\unitlength}{.7mm}
\begin{center}
\begin{picture}(0,0)(150,30) 

\put(48,30){\makebox(0,0){$\scriptstyle (x_1)$}}
\put(48,34){\makebox(0,0){$\scriptstyle 4$}}
\put(51.5,30){\line(1,0){6}}
 \put(57.5,30){\circle*{2}}
\put(58.5,30){\line(1,0){6}}
\put(65.5,30){\circle*{2}}
\put(66.5,30){\line(1,0){6}}
\put(76,30){\makebox(0,0){$\scriptstyle (x_2)$}}
\put(76,34){\makebox(0,0){$\scriptstyle 3$}}
\put(79.5,30){\line(1,0){6}}
\put(86.5,30){\circle*{2}}
\put(87.5,30){\line(1,0){6}}
\put(97,30){\makebox(0,0){$\scriptstyle (x_3)$}}
\put(97,34){\makebox(0,0){$\scriptstyle 2$}}
\put(65.5,29){\line(0,-1){6}}
\put(65.5,22){\circle*{2}}
\put(65.5,21){\line(0,-1){6}}
\put(65.5,16){\circle*{2}}
\put(65.5,15){\line(0,-1){6}}
\put(65.5,07.5){\makebox(0,0){$\scriptstyle (x_4)$}}
\put(71.5,07.5){\makebox(0,0){$\scriptstyle 2$}}

\put(65,0){\makebox(0,0){Figure~1}}
\put(120,22){\makebox(0,0){$\scriptstyle(x_1)$}}
\put(120,26){\makebox(0,0){$\scriptstyle w(x_1)$}}
\put(123.5,22){\line(1,0){6}}
\put(130.5,22){\circle*{2}}
\put(136,22){\makebox(0,0){$\cdots$}}
\put(140,22){\circle*{2}}
\put(141,22){\line(1,0){6}}
\put(150,22){\makebox(0,0){$\scriptstyle(x_2)$}}
\put(150,26){\makebox(0,0){$\scriptstyle w(x_2)$}}
\put(153.5,22){\line(1,0){6}}
\put(160.5,22){\circle*{2}}
\put(166,22){\makebox(0,0){$\cdots$}}
\put(170,22){\circle*{2}}
\put(171,22){\line(1,0){6}}
\put(180,22){\makebox(0,0){$\scriptstyle(x_3)$}}
\put(180,26){\makebox(0,0){$\scriptstyle w(x_3)$}}
\put(183.5,22){\line(1,0){6}}
\put(190.5,22){\circle*{2}}
\put(196,22){\makebox(0,0){$\cdots$}}
\put(202.5,22){\makebox(0,0){$\cdots$}}
\put(207,22){\circle*{2}}
\put(208,22){\line(1,0){6}}
\put(217.5,22){\makebox(0,0){$\scriptstyle(x_n)$}}
\put(217.5,26){\makebox(0,0){$\scriptstyle w(x_n)$}}
\put(166,10){\makebox(0,0){Figure~2}}
\end{picture}

\end{center}
\vspace{3cm}

\begin{example}
\label{ex-cyclic-quotient}{\rm
If $(S,0)$ is a cyclic-quotient singularity i.e. has a graph $\Gamma$ as on figure~2, we may order $\Gamma_{NT}=\{x_1<x_2<\dots<x_n\}$ and each central element $x$ (central vertex  or central arc)  lies in a unique  $[x_i,x_{i+1}]$ and is easily seen to contribute to $\Delta_{S,0}$ by a $\Delta_{x}:=\Delta_{A_{l([x_i,x_{i+1}])}}$, where $l[x_i,x_{i+1}]$ is the number of vertices on the chain $[x_i,x_{i+1}]$ ; the contact between each $\Delta_{x}$ is one (i.e their tangent cones have no common components).
Here $\delta_{x_i}$ is $2 w(x_i)-4$ lines for $i=1$ and $i=n$, and $2 w(x_i)-6$ for $1<i<n$, all this lines being distinct.
So, with $\Delta_{A_n}$ as in notation~\ref{not-delta} i)~: 
$$\Delta_{S,0}=\delta_{x_1}\cup \Delta_{A_{l[x_1,x_2]}}\cup \delta_{x_2}\cup\cdots\Delta_{A_{l[x_{n-1},x_n]}}\cup\delta_{x_n},$$
\noindent with contact one between all the curves in the ``$\cup$''.}
\end{example}

\begin{remark}
\label{rem-BS}{\rm 
In particular, the equisingularity type of $(\Delta_{S,0},0)$ depends only on the resolution graph of $(S,0)$ i.e. of the topological type of $(S,0)$, what is known to be wrong for other normal surface singularities  as shown in \cite{BS}.}
\end{remark}





{\footnotesize
\noindent{\bf Acknowledgement} The origin of this work is an intuition of L\^e's that these discriminants should be describable. I learnt the fine  property~\ref{thm-transv-polaire} from B. Teissier himself. It is a pleasure to thank both of them for all I've learnt from them. I am very grateful to M. Merle and M. Spivakovsky for they pointing out to me respectively the reference \cite{BH} and a mistake in an earlier formulation of this work.}

%

\bigskip
\noindent
C.M.I., 39 rue F. Joliot-Curie,\\
13453 Marseille Cedex 13\\
France\\
email~: bondil@cmi.univ-mrs.fr

\end{document}